\documentclass[a4paper,10pt]{article}
\usepackage{amsmath, amssymb}

\usepackage[francais,english]{babel}

\usepackage{amsmath, amssymb}
\usepackage{geometry,amsmath,amssymb,enumerate,bbm,amsthm}

\usepackage{amsfonts}
\usepackage{amsmath,amsthm,indentfirst}
\usepackage{amssymb}
\usepackage{pstricks}
\usepackage{amsthm}
\usepackage[matrix,arrow]{xy}
\usepackage{amssymb}
\usepackage{amsmath}
\usepackage{vmargin}
\usepackage{dsfont}
\usepackage{mathrsfs}
\usepackage{fancyhdr}

\usepackage[pdftex,bookmarks,colorlinks,breaklinks]{hyperref}

\newcommand{\eps}{{\varepsilon}}
\newcommand{\ep}{\epsilon}

\newcommand{\abs}[1]{\left\vert#1\right\vert}

\newcommand{\dw}{{w\hspace{0,0098889cm}dx}}

\newtheorem{theorem}{Theorem}[section]

\newtheorem{lemma}[theorem]{Lemma}

\theoremstyle{definition}

\newtheorem{remark}[theorem]{Remark}

\numberwithin{equation}{section}

\begin{document}

\date{}
\author{Omar Lazar  
 }
\title{On a 1D nonlocal transport equation with nonlocal velocity and subcritical or supercritical diffusion} 
\maketitle
\bibliographystyle{plain}







\begin{abstract}
 We study a 1D transport equation with nonlocal velocity with subcritical or supercritical dissipation.  For all data in the weighted Sobolev space $H^{k}(w_{\lambda,\kappa}) \cap L^{\infty},$ where  $k=\max(0,3/2-\alpha)$ and $w_{\lambda, \kappa}$ is a given family of Muckenhoupt weights, we prove a global existence result in the subcritical case  $\alpha \in (1,2)$.  We also prove a local existence theorem for large data in $H^{2}(w_{\lambda, \kappa})\cap L^{\infty}$ in the supercritical case $\alpha \in (0,1)$. The proofs are based on the use of the weighted Littlewood-Paley theory, interpolation along with some new commutator estimates.
\end{abstract}

\section{Introduction}
In this paper, we are interested in the following 1D transport equation with nonlocal velocity which was introduced by C\'ordoba, C\'ordoba and Fontelos in \cite{CCF} :
\begin{equation} \label{1d}
\ (\mathcal{T}_{\alpha}) \ : \\\left\{
\aligned
&\partial_{t}\theta+\theta_x \mathcal{H}\theta+ \nu \Lambda^{\alpha}\theta = 0 \hspace{2cm} 
\\ \nonumber
& \theta(0,x)=\theta_{0}(x)
\endaligned
\right.
\end{equation}

\noindent where, $\mathcal{H}$ denotes the Hilbert transform defined by 
$$
\mathcal{H} \theta \equiv \frac{1}{\pi} PV \int \frac{\theta(y)}{y-x } \ dy
$$
and 
\[\Lambda^{\alpha}\theta \equiv (-\Delta)^{\alpha/2} \theta=C_{\alpha}P.V.\int_{\mathbb{R}^{n}}{\frac{\theta(x)-\theta(x-y)}{|y|^{1+\alpha}}dy} \]
where $C_{\alpha}>0$ is a positive constant and $0\leq\alpha< 2$. \\

\noindent This equation is related to a simple scalar model introduced  by  Constantin, Lax and Majda (see \cite{CLM}) in order to get a better understanding of the 3D Euler equation written in terms of the vorticity $w$ (that is the curl of the velocity, and the velocity is determined by $v_x=\mathcal{H}w$), namely the following 1D equation
$$
\partial_t \omega= \omega \mathcal{H} \omega.
$$
In \cite{CLM}, the authors proved that most of the solutions blow-up in finite time. More recently Okamoto, Sakajo and Wunsch \cite{OSW} introduced a generalization of the Constantin, Lax, Majda (CLM) and the De Gregorio's model, namely
$$
\partial_t w + a v w_x - w\mathcal{H}w =0.
$$
The case $a=-1, a=0, a=1$ being respectively the  C\'ordoba, C\'ordoba and Fontelos model (blow-up of regular solutions can occure as shown in \cite{CCF}, \cite{RL}, and \cite{SV}), the CLM model (most solutions blow-up in finite time)  and  the De Gregorio's model (it is conjecured that solutions exist globally in time and numerical evidence presented in \cite{OSW} are in line with the conjecture). This nonlocal transport equation can also be seen as a 1D model of the  dissipative surface quasi geostrophic equation introduced in \cite{CMT} written in a divergence form. It  is therefore closely related to the incompressible 3D Euler equation written in terms of vorticity (see \cite{CMT}). Another motivation is the link between this equation and the Birkoff-Rott equation which modelises the evolution of vortex sheets with surface tension (see e.g.\cite{CCF}). 
As usually, one has to consider three cases depending on the value of $\alpha \in (0,2)$, namely the cases $\alpha<1$, $\alpha=1$, $\alpha > 1$ which are respectively called  supercritical, critical and  subcritical cases.  In the inviscid case ($\eta=0$), C\'ordoba C\'ordoba and Fontelos have shown in \cite{CCF} that there exists a class of smooth initial data for which the solutions blow-up in finite time. This result has been extended to a slightly dissipative case by Li and Rodrigo \cite{RL} who proved the existence of a class of smooth data such that the solutions blow-up in finite time when the dissipation rate is in the range $\alpha \in (0,1/2)$. These results have been proved in four differents way by Silvestre and Vicol  \cite{SV} for $\alpha \in [0,1/2)$ for equation $(\mathcal{T}_{\alpha})$. It is still an open problem to know whether regular solutions blow-up or exist globally in time in the case $1/2\leq \alpha<1$. The eventual regularity of smooth solution in the spirit of previous works for the supercritical SQG equation, has been obtained by Do in \cite{Do}. We also note that the local well posedness in $H^{2}$ has been obtained by Bae and Granero-Belinch\'on in \cite{BG} in the inviscd case. Both the subcritical and critical case are now quite well understood, for instance, estimate are etablished in \cite{Dong}  for data in $H^{3/2-\alpha}$ in the critical case using Littlewood-Paley theory.  Moreover, by adaptading for instance the method  of Constantin and Vicol introduced in \cite{CV} or the approach of Kiselev, Nazarov and Volberg \cite{KNS} one obtain the existence of global smooth solutions in the critical case. Roughly speaking, most of the results known in the critical and subcritical case for the surface quasi-geostrophic equation turn out to work for this 1D model. However, it seems that constructing global  $L^{2}$ solutions in the supercritical or even in the critical case is still not known. Indeed, one can easily derive a nice energy estimate (see  \cite{CCF}) but the lack of compactness prevents one from passing to the weak limit in the nonlinearity and therefore the existence remains an open problem (we refer to \cite{CV} or \cite{LL} for more details). \\

In this article, we study equation $(\mathcal{T}_{\alpha})$ with subcritical or supercritical diffusion. We extended the class of initial data of some  recent existence results (essentially some of those in \cite{Dong}, and \cite{BG}) to the weighted setting. By considering data in weighted Lebesgue or Sobolev spaces with a weight which has the property to be sufficently decaying at infinity say like $\vert x \vert^{-\lambda}$, with $\lambda>0$, one is allowed to consider the problem $(\mathcal{T}_{\alpha})$ with an initial data that behaves for instance like $\theta_{0}(x) \sim \vert x \vert^{-1/2}$ at infinity, and therefore does not necessarily belong to $L^{2}$. In order to avoid integrability problems close to the origin one may consider a weight of the form $w_{\lambda}(x)=(1+\vert x \vert)^{-\lambda}$. Then, since $(\mathcal{T}_{\alpha})$ involves the Hilbert transform, it could be interesting to choose  $w_{\lambda}$ so that the Hilbert transform is a continuous operator in the weighted Lebesgue spaces $L^{p}(w_{\lambda})$, this is equivalent to choose the weight in the Muckenkoupt  class $\mathcal{A}_{p}$ (see \ref{c2} for the definition). Such a weight has not to be integrable far from the origin by definition (it has to obey a reverse H\"older inequality), hence one has to choose $\lambda<1$.
Therefore,  a possible family of weights that one could consider is for instance the one defined by $w_{\lambda}(x)=(1+\vert x \vert)^{-\lambda}$ with  $\lambda \in (0,1)$.  Such kind of weights have been already considered by Farwig and Sohr in \cite{FS}  in the study of the Stokes problem.  More recently, in \cite{LL}, Lemari\'e-Rieusset and the author have studied the critical case for  equation $\mathcal{T}_{\alpha}$ with data in some weighted Sobolev spaces where  the  weight is given by
$w_{\lambda}(x)=(1+x^2)^{-\lambda/2}$.   As a matter of fact, if one would need to control  say $\kappa>1$ derivatives of the weight, then it would be better to use a more general class of radially symmetric weights of the type $w_{\lambda, \kappa}(x)=(1+\vert x \vert^{\kappa})^{-\lambda/\kappa}$ with $\lambda \in (0,1)$ which are still Mukenhoupt weights. In this article,  we would need to control at least two derivatives of the weight thus we shall assume the integer $\kappa\geq2$ to be even (in order to avoid differentiability issues at $-1$).  We shall therefore consider the general family of Muckenhoupt weights given by $w_{\lambda, \kappa}(x)=(1+ \abs{x}^{\kappa})^{-\lambda/\kappa}$.  Note that the weights considered in \cite{FS} and \cite{LL} correspond respectively to the cases $\kappa=1$ and $\kappa=2$.      \\
 
When one try to do {\it{a priori}} estimates in such kind of weighted spaces some difficulties appear. One of the main issue is to control some extra commutators that involve the fractional differential operator $(-\Delta)^{\alpha}$ and the weight $w_{\lambda, \kappa}$. Another obstacle is that the usual Sobolev embedding from $H^{k}$ into  $L^{\infty}$, $k>n/2$  is not any more true in the weighted setting that is why, it must be important to always specify that the data lie in $L^{\infty}$ even if the data is in a sufficiently regular weighted Sobolev space. However, one may still use weighted Sobolev's embedding with intermediate exponents (different from 1 and $\infty$), this is done for instance by Fabes, Kenig and Serapioni in \cite{FKS} or in Maz'ya's book \cite{Maz}.
\newpage
We shall prove a global existence result for data in the weighted Sobolev spaces $H^{k}({w_{\lambda,\kappa}})$ with $k=\max(0,3/2-\alpha)$ for subcritical values of $\alpha$. Beside, we prove a local existence theorem for data in $H^{2}({w_{\lambda,\kappa}})$ in the supercritical case. In both cases, the weight is defined by $w_{\lambda, \kappa}(x)=(1+ \abs{x}^{\kappa})^{-\lambda/\kappa}$ where $\lambda \in (0,1)$ in the subcritical case, and $\lambda \in (0,\alpha/2)$ in the supercritical case, and $\kappa\geq 2$ is an even integer. \\

 A key tool in the proof of those existence results is the use of some new commutator estimates involving the family of weights $w_{\lambda, \kappa}$ and the nonlocal operator $(-\Delta)^{\alpha/2}$ (in particular for supercritical values of $\alpha$). These lemma can be used to treat other equations that involve nonlocal fractional operators. Another tool is the Littlewood-Paley theory in the weighted setting  which allows to deal with data in Sobolev space. It is worth saying that  in the $H^{1}$ or $H^{1/2}$ case, one could follow the approach of \cite{LL} in the subcritical case taking advantage of nice cancellations and formulas involving the Hilbert transform. However, in this paper we aim at stating a theorem that allows to deal with a bigger range of Sobolev regularity. The use of the weighted Littlewood-Paley theory turned out to be efficient to treat the low frequencies of the Hilbert transform which are not continuous in $L^{\infty}$. Also, new commutator estimates are needed to treat the supercritical or subcritical case.  \\

\noindent The article is organized as follows. In the first section, we give the statement of the main theorems. In the third section, we recall some tools and important results that we shall use in the proof of our main theorems. In the third part, we establish {\it{a priori}} estimates for positive initial data  $\theta_0 \in H^{k}( w_{\lambda,\kappa}(x) dx) \cap L^\infty$ with $k=\max(0,3/2-\alpha)$ where the case $k={0}$ and $k\ne0$ are treated separetly into two subsections. The last section is devoted to the proof of the  local existence for $H^{2}_{w}$ data in the supercritical case. 

\section{Main results}

\noindent In the subcritical case, we prove the following global existence result for arbitrary large data in weighted Sobolev spaces.

\begin{theorem} \label{sub}   Assume that $1<\alpha<2$, then for all weights $w_{\lambda,\kappa}(x)=(1+\vert x\vert^{\kappa})^{-\lambda/\kappa}$ with $\lambda \in (0,1)$ and $\kappa\geq2$ an even integer,  and for all positive initial data $\theta_0 \in H^{k}( w_{\lambda, \kappa}(x) dx) \cap L^\infty$ with $k\in \max(0,3/2-\alpha)$, there exists at least one global  solution $\theta$ to the equation $\mathcal{T}_{\alpha},$ which verifies, for all  finite $T>0$
$$
\theta \in \mathcal{C}([0,T], H^{k}(w_{\lambda,\kappa}(x) dx)) \cap L^{2}([0,T] , \dot H^{k+\alpha/2}(w_{\lambda,\kappa}(x) dx)).
$$
Moreover, for all $T<\infty$, we have
$$\Vert \theta (T) \Vert^{2}_{H^{k}(w_{\lambda,\kappa})} \leq \Vert \theta_{0} \Vert^{2}_{H^{k}(w_{\lambda,\kappa})} e^{CT} $$ \\
\vspace{0cm}
\noindent The constant $C>0$ depends  on  $\Vert \theta_{0} \Vert_{L^\infty}$, $\lambda$, $\kappa$ and $\nu$. \\

\end{theorem}

In the supercritical case, we have a local existence result of solutions for arbitrary $H^{2}_{w}$ data.

\begin{theorem} \label{local}   Assume that $0<\alpha<1$, then for all data  $\theta_{0} \in H^{2}_{w_{\lambda,\kappa}}$ where the weight is given by $w_{\lambda}(x)=(1+\vert x\vert^{\kappa})^{-\lambda/\kappa}$ with $\lambda\in (0,\alpha/2)$. Then, there exists a time $T^{*}(\theta_0)>0$ such that  $(\mathcal{T}_{\alpha})$ admits at least one solution that verifies  $$\theta \in \mathcal{C}([0,T], H^{2}(w_{\lambda,\kappa}(x) dx)) \cap L^{2}([0,T], \dot H^{2+\frac{\alpha}{2}}(w_{\lambda,\kappa}(x) dx))$$ for all $T\leq T^{*}$.

\end{theorem}

\begin{remark}
In the proof of Theorem \ref{sub}, we shall treat the case $k=0$ separetly, for the other values of $k$ the proof is based on the use of the weighted Littlewood-Paley  decomposition where a careful treatment of the low frequencies of the velocity is needed. 
\end{remark}

\begin{remark}
Beside extending some previous results from \cite{Dong} and \cite{BG} this article contains also some new commutator estimates that could be of interest in the study of more general nonlocal and nonlinear equations involving the operator $\Lambda^{\alpha}$ that are contained in lemma \ref{lem31}.
\end{remark}

\section{Weighted Littlewood-Paley theory and Muckenhoupt's class \label{intro}}
\vspace{-0.5cm}
 Let $w$ be a positive and locally integrable function. A measurable function $\theta$ is said to belong to the weighted Lebesgue spaces $L^{p}(w)$ (noted also $L^{p}_{w}$ or $L^{p}(w dx)$) with $1\leq p < \infty$ if and only if  
$$
\Vert \theta \Vert_{L^{p}(w)} := \left(\int \vert \theta(x) \vert^{p} \ w(x) \ dx\right)^{1/p} < \infty,
$$
and we have $L^{\infty}=L^{\infty}_{w}$. We shall say that $f$ belongs to the weighted inhomogeneous Sobolev space $H^{s}_{w}$ (or $H_{s} (w)$)  with $\abs{s}<1/2$ if $f \in L^{2}_{w}$ and $\Lambda^{s} \theta \in L^{2}_{w}$, it is endowed with the semi-norm 
$$
\Vert \theta \Vert_{H^{s}_{w}}=\Vert \theta \Vert_{L^{2}_{w}}+ \Vert \Lambda^{s} \theta \Vert_{L^{2}_{w}}. 
$$
Anagolously, we defined the homogeneous weighted Sobolev space $\dot H^{s}_{w}$ (or $ \dot H_{s} (w)$)  as the space such that the following semi-norm is finite
$$
\Vert \theta \Vert_{\dot H^{s}_{w}}= \Vert \Lambda^{s} \theta \Vert_{L^{2}_{w}}. 
$$
We will use the usual notation for the space-time norms, namely, we will say that $\theta$ belongs to the space $L^{2}([0,T], \dot H^{s}(w(x) dx))$ if
$$
\hspace{-1cm}\int_{0}^{T} \int \vert \Lambda^{s}\theta (x,t) \vert^{2} \ w(x) dx \ dt< \infty. 
$$
with the  classical modification for $L^{\infty}([0,T], \dot H^{s}(w(x) dx))$, that is
$$
\hspace{-1cm}\sup_{0<t<T} \int \vert \Lambda^{s} \theta (x,t) \vert^{2} \ w(x)dx < \infty. 
$$

\noindent We shall also use the so-called  Hardy-Littlewood maximal function defined as follows. The  Hardy-Littlewood maximal function of a locally integrable function $\theta$ on  $\mathbb{R}^n$ is defined by
$$\mathcal{M}f (x)=\sup_{Q}\frac{1}{\vert Q \vert} \int_{Q} \vert \theta(y) \vert \ dy,$$
\noindent where the supremum is taken over all cubes $Q$ of  $\mathbb R^n$ centered at the point $x\in \mathbb R^n$  and $\vert Q\vert$ stands for the Lebesgue measure of the cube $Q$. One of the remarkable properties of the maximal Hardy-Littlewood  function is  that it is continuous operator from $L^{p}$ to $L^{p}$ for all $1<p\leq \infty$. For $p=1$, it is a continuous operator  from $L^{1}$ to the weak Lebesgue space $L^{1,\infty}$ endowed with the norm $\Vert \theta \Vert_{L^{1,\infty}}=\sup_{\lambda>0} \left \{{\lambda}^{-1} \vert \{  x \in \mathbb R^n, \vert \theta(x) \vert > \lambda\} \vert \right \}$.

\noindent The extension of this very useful continuity property to weighted Lebesgue spaces goes back to Muckenhoupt \cite{Muc}. He proved that a necessary and sufficient condition on $w$ which ensured the continuity of the  Hardy-Littlewood maximal function on $L^{p}(w)$, $1<p\leq \infty$. More precisely, the Muckenhoupt theorem \cite{Muc} states that there exists a  constant $C_1>0$ such that
\begin{equation} \label{c1}
\quad \quad \quad  \quad \quad \quad \quad \quad \int (\mathcal{M}\theta(x))^{p} w(x) \ dx \leq C_1(w) \int \vert \theta(x) \vert^{p} w(x) \ dx \quad \quad \quad  \quad \quad \quad \quad  
\end{equation}
if and only if, there exists a constant  $C_2 (w)>0$ such that, for all cubes $Q$  in $\mathbb R^{n}$,
\begin{equation}\label{c2}
\quad \quad \quad  \quad \quad \quad \quad \quad \sup_{Q} \left( \frac{1}{\vert Q \vert} \int_Q w \ dx \right) \left( \frac{1}{\vert Q \vert} \int_Q w^{\frac{1}{1-p}} \ dx \right)^{p-1} \leq C_2(w),  \quad \quad \quad  \quad \quad \quad  
\end{equation}

\noindent where $\vert Q\vert$ is the Lebesgue measure of the arbitrary cubes $Q$ with edges parallel to the coordinate axes. Those weights $w$ statisfying  \ref{c2} are said to belong to the $\mathcal{A}_{p} (\mathbb R^n)$ class of Muckenhoupt (\cite{Muc}). The necessary condition is not difficult to obtain, it suffices to set $\theta(x)=w(x)^{\frac{1}{1-p}} \mathds{1}_{Q}(x)$ in \ref{c1}. As for the sufficient condition we refer to \cite{Muc} or \cite {CM2}. \\

\noindent  Another remarkable  property of the $\mathcal{A}_{p}$ weights is  that the Hilbert transform is a continuous operator on the space $L^{p}(w)$, with $1<p<\infty$ if and only if the weight $w$ belongs the $\mathcal{A}_{p}$ class of Muckenhoupt that is those which verify \ref{c2} (see \cite{HMW}). More generally,  this property holds for all Calder\'on-Zygmund operators $T$ (see for instance \cite{CF2}, \cite{CM2}, \cite{St}), namely,  there exists a constant $C(T,w)>0$ which depends on the operator $T$ and on the weight $w$, such that
$$
\Vert T(f) \Vert_{L^{p}(w)} \leq C(T,w)\Vert w \Vert_{\mathcal{A}_p} \Vert f \Vert_{L^{p}(w)}.
$$

\noindent As mentioned in the introduction, the  family of weights $w_{\lambda,\kappa}(x)=(1+\vert x \vert^{\kappa})^{-\lambda/\kappa}$ with  $0<\lambda<1$ and $\kappa\geq2$ being an even integer, that we consider in the article belongs to the $(\mathcal{A}_{p})_{p\in [2,\infty)}$ class of  Muckenhoupt. Therefore, one can use  all the aforementioned useful continuity results.  \\

We finish this section by recalling that, as in the unweighted setting (see for instance \cite{PGLR}, \cite{Can}, \cite{BCD}), the weighted Sobolev $H^{s}_{w}$ spaces can be defined through Littlewood-Paley theory when $w\in \mathcal{A}_{\infty}=\cap_{p>1} \mathcal{A}_{p}$ (see \cite{DSK}). Even in the case $w \in {\mathcal{A}_{p,loc}}$ (that is considering only small cubes in the supremum in \ref{c2}) one still have a satisfactory Littlewood-Paley theory as shown in \cite{Rych}. The construction is as follows, fix a function $\phi_{0} \in \mathcal{D}(\mathbb R^{n})$ that is non-negative and radial and such that ${\phi}_{0}(\xi) = 1$,  if $\vert \xi \vert \leq 1/2$  and ${\phi}_{0}(\xi)=0$  if $\vert \xi \vert \geq 1$. Then, from this fixed function $\phi_{0}$ we define  $\psi_{0}$ so that ${\psi}_{0}(\xi)=\phi_{0}(\xi/2)-\phi_{0}(\xi)$ (which is supported in a corona). For $j \in\mathbb{Z}$, we define the distributions $S_{j}f= \mathcal{F}^{-1}(\phi_{0}(2^{-j}\xi ) \hat{f}(\xi))$ and $\Delta_{j}f= \mathcal{F}^{-1}(\psi_{0}( 2^{-j}\xi) \hat{f}(\xi))$ and we get the so-called (inhomogeneous) Littlewood-Paley decompositon of 
$f \in \mathcal{S}'(\mathbb R^{n})$ that is for all $K\in \mathbb Z$ we have the following inequality in $\mathcal{S}'(\mathbb R^{n})$
\begin{equation} \label{lp}
f=S_{K} f + \sum_{j\geq K} \Delta_{j}f.
\end{equation} 
The homogeneous decomposition is obtained through a passage to the limit in equality \ref{lp} as $K\rightarrow -\infty$ in the $\mathcal{S}'(\mathbb R^{n})$ topology and  we obtain
\begin{equation} \label{lph}
f= \sum_{j\in \mathbb Z} \Delta_{j}f.
\end{equation} 
Equality \ref{lph} is called the homogeneous decomposition and has to be considered modulo polynomials, indeed, $\Delta_{j}f=0 \Leftrightarrow$ $f$ is a polynomial. Then, we define the homogeneous weighted Sobolev spaces $\dot H^{s}_{w}$ for $\vert s \vert <n/2$ as follows
$$
f \in \dot H^{s}_{w}\Longleftrightarrow  f =\sum_{j\in \mathbb Z} \Delta_{j}f,   \ \ {\text {in}} \  \mathcal{S}'(\mathbb R^n)\ \ {\text{and}} \ \  \sum_{j\in \mathbb Z} 2^{2j} \Vert \Delta_{j}f  \Vert^{2}_{L^{2}_{w}}<\infty.
$$
In order to deal with weighted $L^{p}$-estimate of derivatives that involves the operators  $\Delta_{j}$ or $S_{j}$ we shall use the so-called Bernstein's inequality. It says that, for all $f \in \mathcal{S}'(\mathbb{R}^{n})$ and couple  $(j, s) \in \mathbb{Z} \times \mathbb{R}$, and for all $1\leq p \leq q \leq \infty$ and all weights $w \in A_{\infty}$, we have 
 $$\Vert \Lambda^{s} \Delta_{j} f \Vert_{L^{p}_{w}} \lesssim 2^{js} \Vert \Delta_j f \Vert_{L^{p}_{w}} \ \ {\text{also}} \ \ 
 \Vert \Delta_j f \Vert_{L^{q}_{w}} \lesssim2^{{j}(\frac{n}{p}-\frac{n}{q})}  \Vert \Delta_j f \Vert_{L^{p}_{w}} \ \text{and}  \  
\Vert \Lambda^{s} S_j f \Vert_{L^{p}_{w}} \lesssim 2^{js} \Vert S_j f \Vert_{L^{p}_{w}} $$
For two distributions $f$ and $g$ that are in $\mathcal{S}'(\mathbb R^n)$, we may write the paraproduct as follows
$$
fg=\sum_{q \in \mathbb Z} S_{q+1} f \Delta_{q} g + \sum_{j \in \mathbb Z} \Delta_{j} f S_{j} g. 
$$
To establish the existence of at least one solution to equation $(\mathcal{T}_{\alpha})$  it is often useful to truncate the initial data using a function $\psi_{R}$ that is a smooth, positive and compactly supported function in $B_{2R}=[-2R,2R]$. To construct $\psi_R$, we consider a function $\psi \in [0,1]$ such that $\psi(x)=1$ if $\vert x \vert \leq 1$, and $0$ if $\vert x \vert \geq 2$ and then, for $R>0$  we define 
\begin{equation}
\psi_{R}(x)\equiv \psi(x/R).
\end{equation}
This truncation function $\psi_{R}$ will be used throughout the paper when making energy estimates. Note that the constants that will appear in the estimations will always depend on harmless quantities, these constants will be sometimes hidden into the symbol $\lesssim$. Before proving the results, we need to introduce some important lemmas. This is the aim of the next section.

\section{Commutator estimates in the subcritical or supercritical cases.}

In this section, we prove a commutator estimate involving the family of $(\mathcal{A}_p)_{p\in (1,\infty)}$ Muckenhoupt weights $w_{\lambda,\kappa}(x)=(1+\vert x \vert^{\kappa})^{-\lambda/\kappa}$, $0<\lambda<1$, $\kappa\geq2$ being an even integer. More precisely, in our {\it{a priori}} estimates we shall need to control some extra commutators of the type $T_{w_{\lambda,\kappa}}(f)\equiv [\Lambda^{\alpha/2}, w_{\lambda,\kappa}] f$. We shall need to use that such a commutator is continuous operator from $L^{2}(w_{\lambda,\kappa})$ to $ L^{2}(w^{-1}_{\lambda,\kappa})$ in both the subcritical and supercritical cases. In the supercritical case, we need $\lambda$ to be smaller than $\alpha/2$ whereas in the subcritical case we allow $\lambda \in (0,1)$. We shall also need another lemma that gives an $L^{\infty}$ bound for the singular operator $\Lambda^{\alpha}w_{\lambda,\kappa}$ with $\alpha \in (0,2)$. We shall prove both lemmas  in $\mathbb R^n$. The constant $C(n,\alpha, \lambda,\kappa)>0$  that will appear throughout this section can be different from a line to another but we shall keep the same notation for the sake of simplicity.

\begin{lemma} \label{lem31} Let us consider the family of weights given by $w_{\lambda, \kappa}(x)=(1+\abs{x}^{\kappa})^{-\lambda/\kappa}$, then have the following two estimates : \\

$\bullet$ If $\alpha \in (0,1)$, then for all $\lambda$ such that $0<\lambda<\alpha/2$ and for all even integer $\kappa\geq2$, the operator
\begin{eqnarray} \label{comm}
   T_{w_{\lambda, \kappa}} : &L^{2}(w_{\lambda, \kappa})& \longrightarrow L^{2}(w^{-1}_{\lambda, \kappa})  \\
  &f& \nonumber \longmapsto [\Lambda^{\alpha/2}, w_{\lambda, \kappa}] f
\end{eqnarray} 
is continuous. That is, for all $f\in L^{2}(w_{\lambda, \kappa})$, there exists a constant $C>0$ which depends on $\alpha,\lambda, \kappa$, and  $n$ such that,
$$\int \left\vert [\Lambda^{\alpha/2}, w_{\lambda,\kappa}] f \right\vert^{2} \frac{dx}{w_{\lambda, \kappa}} \leq C \int \vert f \vert^{2} w_{\lambda, \kappa} \ dx.$$ 

$\bullet$ If $1<\alpha<2$, then for all $\lambda \in (0,1)$ and for all even integer $\kappa\geq2$, the operator $T_{w_{\lambda, \kappa}}$ is  continuous from $L^{2}(w_{\lambda, \kappa})$ to $L^{2}(w^{-1}_{\lambda, \kappa}).$

\begin{remark}
Note that in the supercritical case, the weight $w_{\lambda, \kappa}$ appearing in the commutator  $T_{w_{\lambda, \kappa}}$  should be denoted  $w_{\lambda(\alpha),\kappa}$ since for each fixed $\alpha$ we have a different weight, but for the sake of simplicity we shall not write the dependence on $\alpha$.
\end{remark}
 
\end{lemma}
\noindent{{\bf Proof of lemma \ref{lem31}}}.  We first remark that
\begin{eqnarray*}
\left\vert \left[\Lambda^{\alpha/2}, w_{\lambda, \kappa}\right] f \right \vert \leq C({\alpha},n) P.V. \int \frac{\vert w_{\lambda, \kappa}(x)-w_{\lambda, \kappa}(y)\vert }{\vert x-y \vert^{n+{\alpha/2}}} \vert f(y) \vert \ dy.
\end{eqnarray*} 
Then, we shall use the following lemma
\begin{lemma} \label{in}
 For all $(x,y) \in \mathbb R^n \times \mathbb R^n $, for all $\lambda \in (0,1)$, and for all even integer $\kappa\geq2$, the following inequality holds
 $$
 \vert w_{\lambda,\kappa}({x})-w_{\lambda,\kappa}({y}) \vert \leq C(\lambda,\kappa) \min\left(\vert {x}-{y}\vert, \vert {x}-{y} \vert^{\lambda/2} \right) \sqrt{w_{\lambda,\kappa}({x})w_{\lambda,\kappa}({y})}.
  $$
\end{lemma}
\noindent{{\bf Proof of lemma \ref{in}.}}  We split the integral into two pieces. If   $x$ and $y$ are such that : $\vert x-y \vert>1$. Then,  if $\vert x-y \vert >\frac{\vert x \vert}{2}$ and $\vert x-y \vert>\frac{\vert y \vert}{2}$ then, since $w_{\lambda,\kappa}({x})-w_{\lambda,\kappa}({y}))\leq w_{\lambda,\kappa}({x}+w_{\lambda,\kappa}({y})$,  it suffices to estimate $w_{\lambda,\kappa}(x)+w_{\lambda, \kappa}(y)$ and in this case
\begin{eqnarray*}
w_{\lambda,\kappa}(x)+w_{\lambda}(y)&\leq& \sqrt{{w_{\lambda, \kappa}(x)}{w_{\lambda,\kappa}(y)}} \left( \frac{1}{\sqrt{w_{\lambda,\kappa}(x)}}+\frac{1}{\sqrt{w_{\lambda, \kappa}(y)}} \right) \\
 &\leq& C(\lambda,\kappa) \vert x-y \vert^{\lambda/2} \sqrt{{w_{\lambda,\kappa}(x)}{w_{\lambda, \kappa}(y)}} 
\end{eqnarray*}
Where we have used that  $(1+2^{\kappa}\vert x - y \vert^{\kappa})^{\lambda/\kappa} > w^{-1}_{\lambda,\kappa}(x) $ then 
$w_{\lambda,\kappa}(x)^{-1/2} \leq (1+2^{\kappa}\vert x-y \vert^{\kappa})^{\lambda/2\kappa} \leq C(\lambda,\kappa) \vert x - y \vert^{{\lambda}/{2\kappa}}$. By symmetry, we also have $w_{\lambda,\kappa}(y)^{-1/2} \leq (1+2^{\kappa}\vert x-y \vert^{\kappa})^{\lambda/2\kappa} \leq C(\lambda, \kappa) \vert x - y \vert^{\lambda/2}$. \\

Otherwise,  $\vert x-y \vert \leq \frac{\vert x \vert}{2}$ or $\vert x-y \vert\leq\frac{\vert y \vert}{2}$. In this case,  we have that $w_{\lambda, \kappa}(x)$ and  $w_{\lambda, \kappa}(y)$ are comparable, in the sense that, one may find a constant $C(\lambda, \kappa)>0$ so that $w_{\lambda, \kappa}(x) \approx C(\lambda,\kappa) w_{\lambda,\kappa}(y)$. Therefore, 
$$
w_{\lambda, \kappa}(x)+w_{\lambda, \kappa}(y) \leq C(\lambda, \kappa) w_{\lambda,\kappa}(x) \leq C(\lambda, \kappa) \sqrt{w_{\lambda, \kappa}(x)w_{\lambda, \kappa}(y)}  \leq C(\lambda, \kappa) \vert x-y\vert^{\lambda/2} \sqrt{w_{\lambda,\kappa}(x)w_{\lambda, \kappa}(y)}
$$
Finally, if $\vert x -y \vert<1$, then 
\begin{eqnarray*}
\vert w_{\lambda, \kappa}(x)-w_{\lambda, \kappa}(y) \vert \leq C \vert x-y \vert \sup_{z \in [x,y]}\vert \nabla w_{\lambda, \kappa}(z) \vert &\leq& C(\lambda, \kappa) \vert x - y \vert \sup_{z\in[x,y]} \vert w_{\lambda, \kappa}(z) \vert \\
 &\leq& C(\lambda, \kappa)\vert x-y \vert \sqrt{w_{\lambda, \kappa}(x)w_{\lambda, \kappa}(y)}.
\end{eqnarray*}
\qed

\noindent Therefore, using lemma \ref{in}, we infer that 
\begin{eqnarray} \label{min}
\left\vert[\Lambda^{\alpha/2}, w_{\lambda, \kappa}] f \right \vert \leq C(\lambda, \kappa,n) \sqrt{w_{\lambda, \kappa}(x)} \int \min \left( \frac{1 }{\vert x-y \vert^{n-1+\frac{\alpha}{2}}}, \frac{1 }{\vert x-y \vert^{{n-\lambda+\frac{\alpha}{2}}}} \right) \sqrt{w_{\lambda, \kappa}(y)}\vert f(y) \vert  \ dy \hspace{0,2cm}
\end{eqnarray}
Let us set 
$$
\chi(x)\equiv\min \left( \frac{1 }{\vert x \vert^{n-1+\frac{\alpha}{2}}}, \frac{1 }{\vert x \vert^{{n-\lambda+\frac{\alpha}{2}}}} \right),
$$
In the supercritical case, we have $0<\alpha<1$, therefore, if $\lambda$ is such that  $0<\lambda<\alpha/2$, then  we observe that  
 $$\chi(x)=\frac{1 }{\vert x \vert^{{n-\lambda+\frac{\alpha}{2}}}} \mathds{1}_{\vert x \vert >1} \in L^{1}(\mathbb R^{n}),$$

\noindent and, 
 $$ \chi(x)=\frac{1 }{\vert x \vert^{n-1+\frac{\alpha}{2}}} \mathds{1}_{\vert x \vert \leq1} \in L^{1}(\mathbb R^{n}).$$
\noindent Hence,  by convolution $x\mapsto \left(\chi * \sqrt{w_{\lambda, \kappa}} f \right)(x) \in L^{2}(\mathbb R^{n})$ and we have for $C=C(\alpha,\lambda, \kappa, n)>0$
\begin{eqnarray} \label{sc}
\left \Vert \int  \min \left( \frac{1 }{\vert x-y \vert^{n-1+\frac{\alpha}{2}}}, \frac{1 }{\vert x-y \vert^{{n-\lambda+\frac{\alpha}{2}}}} \right) \sqrt{w_{\lambda, \kappa}(y)}\ f(y) \ dy \right   \Vert_{L^{2}}  &\leq& C \Vert B \Vert_{L^{1}} \left\Vert \sqrt{w_{\lambda, \kappa}(y)} f(y) \right\Vert_{L^{2}} \nonumber\\
&\leq& C \Vert f \Vert_{L^{2}(w_{\lambda, \kappa})}
\end{eqnarray}
Then, using inequality \ref{min} we conclude that the commutator $f \mapsto T_{w_{\lambda, \kappa}}(f)$  is a continuous operator from $L^{2}(w_{\lambda, \kappa}) \rightarrow L^{2}(w^{-1}_{\lambda, \kappa})$ for supercritical values of $\alpha$.
As for the subcritical case, namely $1<\alpha<2$, we have that, for all $\lambda \in (0,1)$, the function $x \mapsto \chi(x) \in L^{1}(\mathbb R^{n})$. Therefore inequality \ref{sc} still hold, and we may conclude as before.

\qed 

\subsection{Control of the $L^{\infty}$ norm of the nonlocal operator $\Lambda^{\alpha} w_{\lambda}$, with $\alpha \in (0,2)$.} 
\noindent We shall state and prove a general lemma that deals with a bound for $\Lambda^{\alpha} w_{\lambda}$  valid for all values of $\alpha \in (0,2)$. As a matter of fact, we shall only use the estimate for $\alpha=1$. More precisely, we shall only use it in the study of the $L^{2}$ norm, and more precisely when we will integrate by parts and send the $\Lambda$ onto the weight $w_{\lambda, \kappa}$ just after the use of the C\'ordoba and C\'ordoba inequality (see equation \ref{born}). 

\begin{lemma} \label{lem33} For all $\alpha \in (0,2)$,
for all family of weights $w_{\lambda, \kappa}(x)=(1+\vert x \vert^{\kappa})^{-\lambda/\kappa}$, where $0<\lambda<1$ and $\kappa\geq2$ an even integer,  there exists a  constant $C=C(\alpha,\lambda,\kappa, n)>0$ such that the following estimate holds
$$\vert \Lambda^{\alpha} w_{\lambda, \kappa} \vert \leq C  w_{\lambda, \kappa}(x)$$
\end{lemma}
\noindent{{\bf Proof of lemma \ref{lem33}.}} 
It suffices to write
\begin{eqnarray*}
\Lambda^{\alpha} w_{\lambda,\kappa} &&=\frac{1}{2}C(n,\alpha) \int_{\vert x-y \vert \leq 1} 
\frac{2w_{\lambda,\kappa}(x)-w_{\lambda, \kappa}(2x-y)- w_{\lambda, \kappa} (y) }{\vert x - y \vert^{1+\alpha}} \ dy \\
&& \ + \ C(n,\alpha) \int_{\vert x-y \vert>1} \frac{w_{\lambda, \kappa}(x)-w_{\lambda, \kappa}(y)}{\vert x - y \vert^{1+\alpha}} \ dy \\
&&\equiv I_{1}+I_{2}
\end{eqnarray*}
For $I_1$, we see that, using a second order Taylor-expansion 
$$
\vert 2w_{\lambda, \kappa}(x)-w_{\lambda, \kappa}(2x-y)- w_{\lambda, \kappa} (y) \vert \leq \vert x-y \vert^{2}  \sup_{z \in [x,y]}\vert \nabla^{2} w_{\lambda, \kappa}(z) \vert, 
$$
together with the following inequality,
$$
\vert \nabla^{2} w_{\lambda, \kappa}(x) \vert \leq C w_{\lambda,\kappa}(x),
$$
allow us to conclude that,
$$
\vert I_{1}(x) \vert \leq C   w_{\lambda, \kappa}(x).
$$
As for the second integral $I_2$, we split the set  $\Omega(x)=\{y \vert \ \vert x-y \vert>1\}$ into two regions. Namely, we intersect the set $\Omega(x)$ with $$\Omega_1(x) = \{\vert x-y \vert>\frac{\vert x \vert }{2}\}\cap \{\vert x-y \vert > \frac{\vert y \vert}{2} \} \ \ \text{or} \ \ \Omega_2 (x)= \{\vert x-y \vert\leq \frac{\vert x \vert }{2}\}\cup\{\vert x-y \vert \leq \frac{\vert y\vert}{2} \}, $$ and we write
\begin{eqnarray*}
I_{2}&=&C(n,\alpha) \ w_{\lambda, \kappa} (x)\int_{\Omega\cap\Omega_1} \frac{1-\frac{w_{\lambda, \kappa}(y)}{w_{\lambda, \kappa}(x)}}{\vert x - y \vert^{1+\alpha}} \ dy + C(n,\alpha, \kappa) \ w_{\lambda, \kappa}(x)\int_{\Omega\cap\Omega_2} \frac{1-\frac{w_{\lambda, \kappa}(y)}{w_{\lambda, \kappa}(x)}}{\vert x - y \vert^{1+\alpha}} \ dy \\
&\equiv& I_{2,a}+ I_{2,b}.
\end{eqnarray*}
\noindent For the first integral,  we observe that $(1+\vert x \vert^{\kappa})^{\lambda/\kappa} \leq (1+2^{\kappa} \vert x-y \vert^{\kappa})^{\lambda/\kappa} \leq C(\lambda, \kappa) \vert x-y \vert^{\lambda}$ and therefore,  $w^{-1}_{\lambda, \kappa}(x) \leq C(\lambda, \kappa) \vert x-y \vert^{\lambda}$. This allows us to get
\begin{eqnarray}
\nonumber \vert I_{2,a} \vert &\leq&  C  w_{\lambda, \kappa}(x) \left( \int_{\Omega\cap\Omega_1} \frac{1}{\vert x -y \vert^{1+\alpha}}+\int_{\Omega\cap\Omega_1}  \frac{ \vert w^{-1}_{\lambda, \kappa}(x)w_{\lambda, \kappa}(y) \vert}{\vert x-y \vert^{1+\alpha}}  \ dy \right) \\
&\leq& \label{int}C  w_{\lambda, \kappa}(x)\left( C(\alpha) + \int_{\Omega\cap\Omega_1} \frac{1}{\vert x-y \vert^{1+\alpha-\gamma}}  \ dy \right) \\
&\leq&\nonumber C w_{\lambda, \kappa}(x) 
\end{eqnarray}
where we used that $1\leq \alpha$ (critical or subcritical values) or $2\lambda<\alpha<1$ (supercritical value) imply the convergence of the last integral \ref{int} because $0<\lambda<1$.
For the last integral $I_{2,b}$,  we use that $w_{\lambda,\kappa}(x)$ and  $w_{\lambda,\kappa}(y)$ are comparable in the region $\Omega \cap \Omega_{1}$, so that 
$$
\vert I_{2,b} \vert \leq w_{\lambda, \kappa}(x) \ \int_{\Omega\cap\Omega_2} \frac{1}{{\vert x - y \vert^{1+\alpha}}} \ dy  \leq C  w_{\lambda, \kappa}(x)
$$

\qed

 \section{{Global existence for data $\theta_{0} \in H^{k}(w_{\lambda,\kappa}) \cap L^{\infty}$}}
 
 The aim of this section is to prove Theorem \ref{sub}. We shall first treat the case of an $L^{2}_{w} \cap L^{\infty}$ data and then the case of a data in $H^{k}_{w} \cap L^{\infty}$ with $k \ne0$. Indeed, this latter case required a more precise analysis of the low frequencies of the velocity. In the energy estimates, we shall omit to write the dependence of $w$ on $\lambda$ and $\kappa$ for the sake of readibility.

\subsection{The  case  $\theta_{0} \in L^{2}(w_{\lambda,\kappa})\cap L^{\infty}.$  }

 In this subsection, we study the following initial value problem, for all $1<\alpha<2$:

\begin{equation} 
\ (\mathcal{T}_{\alpha}) \ : \\\left\{
\aligned
&\partial_{t}\theta+\theta_x \mathcal{H}\theta+ \nu \Lambda^{\alpha}\theta = 0 \hspace{2cm} 
\\ \nonumber
& 0<\theta(0,x)=\theta_{0}(x)  \in L^{2}( w) \cap L^\infty
\endaligned
\right.
\end{equation}
Before going into the computations to get $\it{a \ priori}$ estimates, one need to check the existence of at least a solution to equation $(\mathcal{T}_{\alpha})$. To do so, we truncate the initial data by multiplying by the function $x\mapsto \phi(x/R)$ that has been introduced at the end of section \ref{intro} and we focus on equation $\mathcal{T}_{{\alpha},R}$  defined as follows
\begin{equation} 
\ (\mathcal{T}_{\alpha,R}) \ : \\\left\{
\aligned
&\partial_{t}\theta_R+\theta_x \mathcal{H}\theta_R+ \nu \Lambda^{\alpha}\theta_R = 0 \hspace{2cm} 
\\ \nonumber
& 0<\theta_{0,R}(x)=\theta_{0}(x) \phi(x/R)  \in L^{2} \cap L^\infty
\endaligned
\right.
\end{equation}

\noindent Let $\theta_{R}$ be a  solution of $\mathcal{T}_{R}$, since the initial data is in $L^2$ and is positive, then it is classical to prove that there exist at least one solution 
$\theta_R \in L^{\infty}([0,T], L^{2}) \cap L^{2}([0,T], \dot H^{\alpha/2})$ (see \cite{CCF}). Then multiplying equation  $\mathcal{T}_{R}$ by such a solution $\theta_R$ and integrating in space give
$$
\partial_{t} \int \frac{\theta^{2}_{R}}{2}w \ dx = -\int \theta_R \partial_{x} \theta_R \mathcal{H}\theta_{R} w \ dx -  \int  w \theta_{R} \Lambda^{\alpha} \theta_{R} \ dx,
$$
which can be rewritten as,
$$
\partial_{t} \int \frac{\theta^{2}_{R}}{2}w \ dx = -\int \theta_{R} \partial_{x} \theta_{R} \mathcal{H}\theta_{R} w \ dx -  \int \Lambda^{\alpha/2} \theta_{R} [ \Lambda^{\alpha/2}, w] \theta_{R} \ dx  - \int \vert \Lambda^{\alpha/2} \theta_{R} \vert^{2} w \ dx.
$$
An integration by parts, together with the fact that $\partial_x \mathcal{H}=- \Lambda$,  give
$$
\partial_{t} \int \frac{\theta^{2}_{R}}{2}w \ dx = - \int \theta^{2}_{R} \Lambda \theta_{R} w \ dx \ + \int \theta_{R} \mathcal{H}\theta w_x \ dx \ -  \int \Lambda^{\alpha/2} \theta_{R} [w, \Lambda^{\alpha/2}] \theta_{R} \ dx \ -\int \vert \Lambda^{\alpha/2} \theta_{R} \vert^{2} w \ dx.
$$
By using the  C\'ordoba and C\'ordoba inequality \cite{CC} : $\Lambda \theta^{3}_{R} \leq 3\theta^{2}_{R} \Lambda \theta_{R}$ \ (the convexity of $\theta \mapsto \theta^{3}_{R}$ is ensured by the positivity assumption on $\theta_{R}$), and then integrating by parts yields
\begin{eqnarray}
 \label{born}\partial_{t} \int \frac{\theta^{2}_{R}}{2}w \ dx \leq - \frac{1}{3} \int   \theta^{3}_{R} \Lambda w \ dx  + \int \theta_{R} \mathcal{H}\theta_{R} w_x \ dx &-&  \int \Lambda^{\alpha/2} \theta_{R} [w, \Lambda^{\alpha/2}] \theta_{R} \ dx \\
&-& \nonumber\int \vert \Lambda^{\alpha/2} \theta_{R} \vert^{2} w \ dx.
\end{eqnarray}
Using lemma \ref{lem33} and the fact that pointwisely  we have the inequality $\partial_{x} w< w$ (indeed,  it suffices to see that  $\partial_{x}w_{\lambda,\kappa}(x)=  -\lambda  x^{k-1} (1+\vert x \vert^{\kappa})^{-1}
w_{\lambda,\kappa}(x)$ and since $0<\lambda<1$, we obviously have ${-{\lambda}} x^{k-1} < 1+\vert x \vert^{\kappa}$), we obtain
\begin{eqnarray*}
 \partial_{t} \int \frac{\theta^{2}_{R}}{2}w \ dx&\leq& \frac{1}{3}\int \vert \theta_{R} \vert \vert \theta_{R} \sqrt{w} \vert^{2} \ dx \ + \int \vert \sqrt{w}\theta_{R} \vert \vert  \mathcal{H}\theta_{R} \ \sqrt{w} \vert \ dx \ + \int \vert \sqrt{w}\Lambda^{\alpha/2} \theta_{R} \vert \vert [\Lambda^{\alpha/2},w ] \theta_{R} \vert \frac{dx}{\sqrt{w}} \\ &&\ -\int \vert \Lambda^{\alpha/2} \theta_{R} \vert^{2} w \ dx.
\end{eqnarray*}
Recalling that $\Vert \theta_{R}\Vert_{L^{2}(w)}=\Vert \theta_{R} \sqrt{w} \Vert_{2}$ and $\Vert \theta_{R}\Vert_{\dot H^{\alpha/2}(w)}=\Vert  \sqrt{w}  \Lambda^{\alpha/2} \theta_{R}\Vert_{2}$, and using lemma \ref{lem31} together with the $L^\infty$ maximum principle,  and the continuity of the Hilbert transform on $L^{2}(w)$ (because $w \in \mathcal{A}_{2}$) we finally get
$$
\partial_{t} \int \frac{\theta^{2}_{R}}{2}w \ dx \leq C\Vert \theta_{0,R} \Vert_{L^\infty} \Vert \theta_{R} \Vert^{2}_{L(w)} +\Vert \Lambda^{\alpha/2} \theta_{R} \Vert_{L^{2}(w)} \Vert \theta_{R} \Vert_{L^{2}(w)}- \Vert \Lambda^{\alpha/2}\theta_{R} \Vert^{2}_{L^{2}(w)}.
$$
Using Young's inequality we infer that, for all $\eta>0$
\begin{equation}
\partial_{t} \int \frac{\theta^{2}_{R}}{2}w \ dx \leq \left(\frac{2}{\eta}+C\Vert \theta_{0,R} \Vert_{L^\infty}\right) \Vert \theta \Vert^{2}_{L^{2}(w)}+\left(\frac{\eta}{2}-1\right)  \Vert \Lambda^{\alpha/2} \theta_{R} \Vert^{2}_{L^{2}(w)}.
\end{equation}
By choosing $\eta$ sufficiently small (for instance $\eta<2$), we find
\begin{equation*}
\partial_{t} \int \frac{\theta^{2}_{R}}{2}w + C \Vert \Lambda^{\alpha/2} \theta_{R} \Vert^{2}_{L^{2}(w)} dx \leq C\Vert \theta_{0} \Vert_{L^\infty} \Vert \theta_{R} \Vert^{2}_{L^{2}(w)}.
\end{equation*}
Integrating in time $s\in [0,T]$ yields
$$
\Vert \theta_R (x, T) \Vert^{2}_{L^{2}(w)} + C \int_{0}^{T} \Vert \Lambda^{\alpha/2} \theta_{R} \Vert^{2}_{L^{2}(w)} \ ds \leq \Vert \theta_{0,R} \Vert^{2}_{L^{2}(w)} +   C\Vert \theta_{0,R} \Vert_{L^\infty} \int_{0}^{T} \Vert \theta_{R}(x,s) \Vert^{2}_{L^{2}(w)} \ ds.
$$
Hence, by Gronwall's inequality we get 

$$
\hspace{-1cm}\sup_{0<t<T} \int \theta^{2}_{R} w \ dx < \infty \ \ \text{and} \ \ \int_{0}^{T} \int \vert \Lambda^{\alpha/2} \theta_{R} \vert^{2} w \ dx < \infty.
$$
In particular, we have, for all $T<\infty$
$$
\Vert \theta_{R} (T) \Vert^{2}_{L^{2}(w)} \leq \Vert \theta_{0,R} \Vert^{2}_{L^{2}(w)}  e^{C T},
$$ 
where $C>0$  depends only on $\Vert \theta_{0,R} \Vert_{L^{\infty}}$, and $\beta$.  Then, the passage to the weak limit as $R\rightarrow \infty$ allows us to prove the first statement of Theorem \ref{sub}.  Indeed, since $1<\alpha<2$, one has enough compactness to pass to the weak limit in the nonlinear term. This can be done  for instance via Rellich's compactness theorem \cite{PGLR} which provides the strong convergence in $(L^{2}L^{2})_{loc}$ of $\theta_R$ and allows us to pass to the weak limit in the nonlinear term. It is obvious that $\theta_{0,R}$ converges strongly in $L^{2}$ as well as in $L^{\infty}$ toward $\theta_{0}$. This concludes the proof. 

\qed

It is worth recalling that this is an open problem in the critical or super-critical case and when the data is just in $L^{2}$ or $L^{2}_{w}$. Although one can derive nice {\it{a\ priori}} estimates, it is still not clear whether one can  to pass to the limit because of the lack of compactness (see for instance \cite{SV}, \cite{LL}).

\subsection{The  case $k=\max(0,3/2-\alpha)$.}

\noindent In this section, we study the following initial value problem 

\begin{equation} 
\ (\mathcal{T}_{\alpha}) \ : \\\left\{
\aligned
&\partial_{t}\theta+\theta_x \mathcal{H}\theta+ \nu \Lambda^{\alpha}\theta = 0 \hspace{2cm} 
\\ \nonumber
& 0<\theta(0,x)=\theta_{0}(x) \in H^{k}(w) \cap L^\infty. 
\endaligned
\right.
\end{equation}
We shall focus on the following approximate family of equations (where $\phi$ have been introduced in section \ref{intro}).
\begin{equation} 
\ (\mathcal{T}_{\alpha,R}) \ : \\\left\{
\aligned
&\partial_{t}\theta_R+\partial_{x}\theta_R \mathcal{H}\theta_R+ \nu \Lambda^{\alpha}\theta_R = 0 \hspace{2cm} 
\\ \nonumber
& 0<\theta_{0,R}(x)=\theta_{0}(x)\phi(x/R) \in H^{k} \cap L^\infty.
\endaligned
\right.
\end{equation}
From theorem 3.1 of \cite{CCF}, we know that there exists at least one  solution $\theta_{R} \in L^{\infty}L^{2} \cap L^{2}H^{k+\alpha/2}$. For the sake of readibility, we shall omit to write the index $R>0$ in the energy estimates. In order to prove our result, we shall consider the evolution of the following weighted Sobolev norm
$$
\Vert \theta \Vert^{2}_{ H^{k}(w)} =  \int \vert \Lambda^{k} \theta \vert^{2} \ w(x) dx \ + \ \int \vert \theta \vert^{2} \ w(x) dx  =  \Vert \sqrt{w} \Lambda^{k} \theta \Vert^{2}_{L^{2}} + \Vert \sqrt{w}  \theta \Vert^{2}_{L^{2}}.
$$
The main goal is to show that, for all finite $T$, we have
$$
 \Vert \theta(T)\Vert^{2}_{H^{k}(w)} + C \int_{0}^{T} \int \vert \Lambda^{\frac{\alpha}{2}+k} \theta \vert^2 w \ dx \ ds \lesssim \Vert \theta_0\Vert^{2}_{H^{k}(w)} + \int_{0}^{T} \Vert \theta(x,s) \Vert^{2}_{H^{k}(w)} \ ds.
$$
In the previous section, we have already seen that, for all $\eta>0$
\begin{equation} \label{L2b}
\partial_{t} \int \frac{\theta^{2}}{2}w \ dx \leq \left(\frac{2}{\eta}+C\Vert \theta_{0} \Vert_{L^\infty}\right) \Vert \theta \Vert^{2}_{L^{2}(w)}+\left(\frac{\eta}{2}-1\right)  \Vert \Lambda^{\alpha/2} \theta \Vert^{2}_{L^{2}(w)}
\end{equation}
Therefore, it suffices to  focus on the evolution of the homogeneous Sobolev norm. One can write it in terms of controlled commutators (via lemma \ref{lem31}) as follows 
\begin{eqnarray*}
\frac{1}{2} \partial_{t} \int  \vert \Lambda^{k} \theta \vert^{2} \ w dx &=& -\int  \Lambda^{k} \theta  \Lambda^{k}(\theta_{x} \mathcal{H} \theta) w \ dx  -\nu \int \Lambda^{\frac{\alpha}{2}} (w\Lambda^{1/2} \theta) \ \Lambda^{\frac{\alpha}{2} + k } \theta \ dx \\
&=&  -\int \Lambda^{\sigma}( w \Lambda^{k} \theta)  \Lambda^{k-\sigma}(\theta_{x} \mathcal{H} \theta) \ dx  - \nu \int \sqrt{w} \Lambda^{\frac{\alpha}{2} + k  } \theta \  [\Lambda^{\frac{\alpha}{2}}, w] \Lambda^{k}\theta \ \frac{dx}{\sqrt{w}} \\ && \hspace{4cm}   - \nu \int \vert \Lambda^{\frac{\alpha}{2} + k }\theta \vert^{2} \ w dx \\
&=& -\int \sqrt{w} \Lambda^{k-\sigma}(\theta_{x} \mathcal{H} \theta) \frac{1}{\sqrt{w}}[\Lambda^{\sigma}, w] \Lambda^{k} \theta \ dx - \int \sqrt{w} \Lambda^{k+\sigma} \theta \sqrt{w}   \Lambda^{k-\sigma}(\theta_{x} \mathcal{H} \theta) \ dx \\
&&- \nu \int \sqrt{w} \Lambda^{\frac{\alpha}{2} + k  } \theta \  [\Lambda^{\frac{\alpha}{2}}, w] \Lambda^{k}\theta \ \frac{dx}{\sqrt{w}} - \nu \int \vert \Lambda^{\frac{\alpha}{2} + k }\theta \vert^{2} \ w dx \\
&=& (I) + (II) + (III)+(IV)
\end{eqnarray*}

\noindent To estimate $(I)$ and $(II)$ we use the dyadic Littlewood-Paley decomposition in the weighted setting. We have, using the paraproduct formula 
$$
\Vert \theta_{x} \mathcal{H} \theta \Vert_{\dot H^{k-\sigma}(w)} \leq \sum_{q \in \mathbb Z} 2^{q(k-\sigma)} \Vert S_{q+1}\theta_{x} \Delta_{q} \mathcal{H}\theta \Vert_{L^{2}(w)} +  \sum_{j \in \mathbb Z} 2^{j(k-\sigma)} \Vert \Delta_{j}\theta_{x} S_{j} \mathcal{H}\theta \Vert_{L^{2}(w)}.
$$
To estimate the first sum, we use  Bernstein's inequality and  the continuity of the Hilbert transform on $L^{2}(w)$ because $w\in \mathcal{A}_2$, we get
\begin{eqnarray*}
 \sum_{q\in \mathbb{Z}} 2^{q(k-\sigma)} \Vert S_{q+1}\theta_{x} \Delta_{q} \mathcal{H}\theta \Vert_{L^{2}(w)} &\leq&   \sum_{q \in \mathbb{Z}} 2^{q(k-\sigma)} \Vert S_{q+1}\theta_{x} \Vert_{L^{\infty}}  \Vert \Delta_{q} \mathcal{H}\theta \Vert_{L^{2}(w)} \\
 &\leq&  \Vert \theta\Vert_{L^{\infty}}   \sum_{q \in \mathbb{Z}} 2^{(k+1-\sigma)q}   \Vert \Delta_{q} \mathcal \theta \Vert_{L^{2}(w)}  \\
  &\leq&  \Vert \theta_{0}\Vert_{L^{\infty}}   \Vert \theta \Vert_{\dot H^{k+1-\sigma}(w)}.
\end{eqnarray*}
As for the second sum, we need to interpolate since we do not control the $L^{\infty}$ norm of the low frequencies of $\mathcal{H}\theta$. More precisely, we use the fact that $\theta \in L^{2}(w) \cap L^{\infty}$ then by interpolation $\mathcal{H}\theta \in L^{s}(w)$  for all $2<s<\infty$ and we have
$$
\Vert \mathcal{H} \theta \Vert_{L^{s}(w)} \leq C \Vert \theta \Vert^{1-2/s}_{L^{\infty}} \Vert \theta \Vert^{2/s}_{L^{2}(w)}.
$$
Then, if $r$ is the real such that 
$$
\frac{1}{s} + \frac{1}{r}= \frac{1}{2},
$$
\noindent we obtain, using H\"older and then Bernstein's inequality, that
\begin{eqnarray*}
  \sum_{j\in \mathbb{Z}} 2^{j(k-\sigma)} \Vert \Delta_{j}\theta_{x} S_{j} \mathcal{H}\theta \Vert_{L^{2}(w)} &\leq&   \sum_{j\in \mathbb{Z}} 2^{j(k-\sigma)}  \Vert \Delta_{j}\theta_x \Vert_{L^{r}(w)} \Vert S_{j} \mathcal{H}\theta \Vert_{L^{s}(w)} \\
  &\leq&  C   \sum_{j\in \mathbb{Z}} 2^{j(k-\sigma)}   2^{j(\frac{1}{2}-\frac{1}{r})} 2^{j} \Vert \Delta_{j}\theta \Vert_{L^{2}(w)} \Vert  \mathcal{H}\theta \Vert_{L^{s}(w)}  \\
    &\leq&  C  \Vert \theta \Vert^{1-2/s}_{L^{\infty}} \Vert \theta \Vert^{2/s}_{L^{2}(w)}  \sum_{j\in \mathbb{Z}} 2^{j({k-\sigma+\frac{1}{s}+1})} \Vert \Delta_{j}\theta \Vert_{L^{2}(w)}  \\
  &\leq& C \Vert \theta_{0} \Vert^{1-2/s}_{L^{\infty}} \Vert \theta \Vert^{2/s}_{L^{2}(w)} \Vert \theta \Vert_{\dot H^{k+\frac{1}{s}+1-\sigma}(w)}.
\end{eqnarray*}
Therefore,
\begin{eqnarray*}
\Vert \theta_{x} \mathcal{H} \theta \Vert_{\dot H^{k-\sigma}(w)} \leq \Vert \theta_{0}\Vert_{L^{\infty}}   \Vert \theta \Vert_{\dot H^{k + 1-\sigma}(w)} +   \Vert \theta_{0} \Vert^{1-2/s}_{L^{\infty}} \Vert \theta \Vert^{2/s}_{L^{2}(w)} \Vert \theta \Vert_{\dot H^{k+\frac{1}{s}+1-\sigma}(w)}.
\end{eqnarray*}
Then,  lemma \ref{lem31} gives 
$$
(I) \lesssim \Vert \theta_{0}\Vert_{L^{\infty}} \Vert \theta \Vert_{\dot H^{k}(w)}      \Vert \theta \Vert_{\dot H^{k +1-\sigma}(w)} +    \Vert \theta_{0} \Vert^{1-2/s}_{L^{\infty}} \Vert \theta \Vert_{\dot H^{k}(w)} \Vert \theta \Vert^{2/s}_{L^{2}(w)} \Vert \theta \Vert_{\dot H^{k+\frac{1}{s}+1-\sigma}(w)}.
$$
Choosing $s$ and $\sigma$  such that $\frac{1}{s}+1-\sigma< \frac{\alpha}{2}$   (note that it suffices to choose $s$ sufficiently big and $\sigma>\frac{1}{2}$, actually we will choose $\sigma=\frac{1}{2}+\frac{1}{s}$), and using the following interpolation inequality valid for all $\gamma \in (0,1)$ $$\Vert \theta \Vert_{\dot H^{k+\frac{1}{s}+1-\sigma}(w)} \leq \Vert \theta \Vert^{1-\gamma}_{\dot H^{k+\frac{1}{s}+1-\sigma-\gamma}(w)}\Vert \theta \Vert^{\gamma}_{\dot H^{k+\frac{1}{s}+2-\sigma-\gamma}(w)},$$ we obtain
\begin{eqnarray*} 
(I) &\leq& \Vert \theta_{0}\Vert_{L^{\infty}} \Vert \theta \Vert_{\dot H^{k}(w)}      \Vert \theta \Vert_{\dot H^{k + \frac{\alpha}{2}-\frac{1}{s}}(w)} +    \Vert \theta_{0} \Vert^{1-2/s}_{L^{\infty}} \Vert \theta \Vert_{\dot H^{k}(w)} \Vert \theta \Vert^{2/s}_{L^{2}(w)} \Vert \theta \Vert_{\dot H^{k+\frac{1}{s}+1-\sigma}(w)} \\
&\lesssim&  \Vert \theta \Vert_{ H^{k}(w)}      \Vert \theta \Vert_{ H^{k + \frac{\alpha}{2}-\frac{1}{s}}(w)} +   \Vert \theta \Vert^{2/s}_{L^{2}(w)}  \Vert \theta \Vert_{ H^{k}(w)}  \Vert \theta \Vert^{1-\gamma}_{ H^{k+\frac{1}{s}+1-\sigma-\gamma}(w)}\Vert \theta \Vert^{\gamma}_{ H^{k+\frac{1}{s}+2-\sigma-\gamma}(w)}.  
\end{eqnarray*}
Then, we use that $ H^{k}({w}) \hookrightarrow   H^{k+\frac{1}{s}+1-\sigma-\gamma}({w})$ (note that it suffices to choose $\gamma=1-\eps$ with $\ep$ small enough, so that $\ep + \frac{1}{s} < \sigma$) 
moreover we have  $  H^{k+\frac{1}{s}+2-\sigma-\gamma}(w)  \hookrightarrow  H^{k+ \frac{\alpha}{2}}(w)$ (since $\sigma=\frac{1}{2}+\frac{1}{s}$ and we choose $s$ big enough and $\ep$ small enough so that  $\frac{1}{2}+\ep < \frac{\alpha}{2}$ therefore $\frac{1}{s}+ \ep+1-\sigma < \frac{\alpha}{2}$). We get
\begin{equation} {\label{eq: Inegal}}
(I) \lesssim  \Vert \theta \Vert_{ H^{k}(w)}      \Vert \theta \Vert_{ H^{k + \frac{\alpha}{2}-\frac{1}{s}}(w)} +   \Vert \theta \Vert^{2/s}_{L^{2}(w)}    \Vert \theta \Vert^{2-\gamma}_{ H^{k}(w)}\Vert \theta \Vert^{\gamma}_{ H^{k+\frac{\alpha}{2}}(w)}.
\end{equation}

\noindent Then, we shall repeatedly use  Young's inequality, first with the exponent $p_1=\frac{1}{1-\gamma}>1$,  we obtain 
\begin{eqnarray*}
(I)  &\lesssim&     \Vert \theta \Vert_{ H^{k}(w)}      \Vert \theta \Vert_{ H^{k + \frac{\alpha}{2}-\frac{1}{s}}({w})}  +  (1-\gamma)\Vert \theta \Vert^{\frac{1}{s(1-\gamma)}}_{L^{2}({w})}    \Vert \theta \Vert_{ H^{k}({w})} + \gamma \Vert \theta \Vert_{ H^{k+\frac{\alpha}{2}}(w)} \Vert \theta \Vert^{\frac{1}{s\gamma}}_{L^{2}({w})}   \Vert \theta \Vert^{\frac{1}{\gamma}}_{H^{k}(w)}, 
\end{eqnarray*}
and then in the 3 terms of the above inequality (respectively, with the exponents  $p_2= 2$ for the first two ones, and  $p_3=1+\gamma>1$ for the last term), we find that for all $\mu_1>0$
\begin{eqnarray*}
(I) &\lesssim&  \frac{1}{2\mu_1}  \Vert \theta \Vert^{2}_{ H^{k}(w)} +  \frac{\mu_1}{2} \Vert  \theta \Vert^{2}_{ H^{k+\frac{\alpha}{2}}(w)} +  \frac{(1-\gamma)}{2}\Vert \theta \Vert^{\frac{2}{s(1-\gamma)}}_{L^{2}(w)} +  \frac{(1-\gamma)}{2 } \Vert \theta \Vert^{2}_{ H^{k}(w)}   \\
 &&  \ + \ \frac{\mu^{1+\gamma}_{2}}{(1+\gamma) }  \Vert \theta \Vert^{1+\gamma}_{ H^{k+\frac{\alpha}{2}}(w)} +  \frac{\gamma}{(1+\gamma)\mu^{1+\frac{1}{\gamma}}_{2}} \Vert \theta \Vert^{\frac{1}{s(1+\gamma)}}_{L^{2}(w)} \Vert \theta \Vert^{\frac{1}{1+\gamma}}_{ H^{k}(w)}.   
\end{eqnarray*}
Then,  using once again Young's inequality in the last term of the previous inequality (with the exponent $p_4= 2+2\gamma>1$) we finally get
\begin{eqnarray*}
(I)&\lesssim&   \frac{1}{2\mu_1}  \Vert \theta \Vert^{2}_{\dot H^{k}(w)} +  \frac{\mu_1}{2} \Vert  \theta \Vert^{2}_{ H^{k+\frac{\alpha}{2}}(w)}  +  \frac{(1-\gamma)}{2}\Vert \theta \Vert^{\frac{2}{s(1-\gamma)}}_{L^{2}(w)} +  \frac{(1-\gamma)}{2} \Vert \theta \Vert^{2}_{ H^{k}(w)}   \\
 &&  \ + \ \frac{\mu^{1+\gamma}_{2}}{(1+\gamma) } \Vert \theta \Vert^{1+\gamma}_{ H^{k+\frac{\alpha}{2}}(w)} +   \frac{\gamma(1+2\gamma)}{2\mu_{2}^{1+\frac{1}{\gamma}}(1+\gamma)^2}    \Vert \theta \Vert^{\frac{2}{s(1+2\gamma)}}_{L^{2}(w)} + \frac{\gamma}{2\mu_{2}^{1+\frac{1}{\gamma}}(1+\gamma)^2}  \mu_{2}^{1+\frac{1}{\gamma}} \Vert \theta \Vert^{2}_{ H^{k}(w)}.
\end{eqnarray*}
Since $s$ is chosen big enough, we have $\frac{2}{s(1+2\gamma)}<2$ and $\frac{2}{s(1-\gamma)}<2$ (note that the first inequality holds for instance if $s>1$, whereas the second is verified if $s >\frac{2}{\ep}$). Actually, 
in the estimation of $(II)$ below we will need to assume that  $\sigma=\frac{1}{s} + \frac{1}{2}\leq\frac{\alpha}{2}$; thus all those conditions imply that  $s>\max(\frac{2}{\ep}, \frac{2}{\alpha-1})$. Futhermore, we obviously have that $1+\gamma<2$,  hence we obtain, 

\begin{eqnarray*}
(I) &\lesssim&   \frac{1}{2\mu_1}  \Vert \theta \Vert^{2}_{ H^{k}(w)} +  \frac{\mu_1}{2} \Vert  \theta \Vert^{2}_{ H^{k+\frac{\alpha}{2}}(w)}  +  \frac{(1-\gamma)}{2}\Vert \theta \Vert^{2}_{L^{2}(w)} +  \frac{(1-\gamma)}{2} \Vert \theta \Vert^{2}_{ H^{k}(w)}   \\
 &&  \ + \ \frac{\mu^{1+\gamma}_{2}}{(1+\gamma)}  \Vert \theta \Vert^{2}_{ H^{k+\frac{\alpha}{2}}(w)} +   \frac{\gamma(1+2\gamma)}{2 \mu_{2}^{1+\frac{1}{\gamma}} (1+\gamma)^2}\Vert \theta \Vert^{2}_{L^{2}(w)} + \frac{\gamma}{2 \mu_{2}^{1+\frac{1}{\gamma}}(1+\gamma)^2}\Vert \theta \Vert^{2}_{ H^{k}(w)}.
\end{eqnarray*} 

\noindent Hence, we have obtained
\begin{eqnarray*}
(I)&\lesssim&  \left(\frac{\mu_1}{2}  + \frac{\mu^{1+\gamma}_{2}}{1+\gamma} \right) \Vert  \theta \Vert^{2}_{ H^{k+\frac{\alpha}{2}}({w)}} +  \left(  \frac{1}{2\mu_1}  + \frac{1-\gamma}{2}  +\frac{\gamma}{2 \mu_{2}^{1+\frac{1}{\gamma}} (1+\gamma)^2} \right) \Vert \theta \Vert^{2}_{ H^{k}(w)} \\  &+&  \left( \frac{1-\gamma}{2}  +   \frac{\gamma(1+2\gamma)}{2 \mu_{2}^{1+\frac{1}{\gamma}} (1+\gamma)^2}  \right) \Vert \theta \Vert^{2}_{L^{2}(w)}.
\end{eqnarray*}
For $(II)$, we have
$$
(II) \leq \Vert \theta_{0}\Vert_{L^{\infty}}  \Vert \theta \Vert_{\dot H^{k+\sigma}(w)}  \Vert \theta \Vert_{\dot H^{k +1-\sigma}(w)} +  \Vert \theta_{0} \Vert^{1-2/s}_{L^{\infty}} \Vert \theta \Vert_{\dot H^{k+\sigma}(w)} \Vert \theta \Vert^{2/s}_{L^{2}(w)} \Vert \theta \Vert_{\dot H^{k+\frac{1}{s}+1-\sigma}(w)}.
$$ 
Therefore, using the fact that $s$ is chosen so that $\sigma=\frac{1}{s} + \frac{1}{2}\leq \frac{\alpha}{2}$ and  that $1-\gamma=\ep$, we get 
\begin{eqnarray*}
&(II)& \lesssim  \Vert \theta \Vert_{\dot H^{k+\sigma}(w)}  \Vert \theta \Vert_{\dot H^{k + 1- \sigma }(w)} +  \Vert \theta \Vert_{\dot H^{k+\sigma}(w)} \Vert \theta \Vert^{2/s}_{L^{2}(w)}  \Vert \theta \Vert^{1-\gamma}_{\dot H^{k+\frac{1}{s}+1-\sigma-\gamma}(w)}\Vert \theta \Vert^{\gamma}_{\dot H^{k+\frac{1}{s}+2-\sigma-\gamma}(w)} \\
&& \lesssim  \Vert \theta \Vert_{ H^{k+\frac{\alpha}{2}}(w)}  \Vert \theta \Vert_{ H^{k + \frac{1}{2}-\frac{1}{s}}(w)} +  \Vert \theta \Vert_{ H^{k+\frac{\alpha}{2}}(w)} \Vert \theta \Vert^{2/s}_{L^{2}(w)}  \Vert \theta \Vert^{1-\gamma}_{ H^{k+\frac{1}{s}+1-\sigma-\gamma}(w)}\Vert \theta \Vert^{\gamma}_{ H^{k+\frac{1}{s}+2-\sigma-\gamma}(w)}.
\end{eqnarray*}
Then, using the interpolation inequality
$$
\Vert \theta \Vert_{ H^{k+\frac{1}{2}-\frac{1}{s}}(w)} \leq \Vert \theta \Vert^{1-\frac{1}{s}}_{ H^{k+\frac{1}{2}}(w)}  \Vert \theta \Vert^{\frac{1}{s}}_{ H^{k-\frac{1}{2}}(w)},
$$
we obtain,
\begin{eqnarray*}
(II)\lesssim \Vert \theta \Vert^{2-\frac{1}{s}}_{ H^{k+\frac{\alpha}{2}}}  \Vert \theta \Vert^{1/s}_{ H^{k}(w)} +  \Vert \theta \Vert_{ H^{k+\frac{\alpha}{2}}(w)} \Vert \theta \Vert^{2/s}_{L^{2}(w)}  \Vert \theta \Vert^{1-\gamma}_{
 H^{k+\frac{1}{s}+1-\sigma-\gamma}(w)}\Vert \theta \Vert^{\gamma}_{ H^{k+\frac{1}{s}+2-\sigma-\gamma}(w)}.
\end{eqnarray*}
Then, we use Young's inequality in the first term of the above inequality (with the exponent $p_4=2s>1$) and we estimate the last term as we did for $(I)$ (in particular, we use the same embeddings), we get, for all $\mu_3>0$, 
\begin{eqnarray*}
&(II)&\lesssim \frac{2s-1}{2s} {\mu_3^{\frac{2s}{2s-1}}}  \Vert \theta \Vert^{2}_{ H^{k+\frac{\alpha}{2}}(w)} + \frac{1}{2s\mu_{3}^{2s}}   \Vert \theta \Vert^{2}_{ H^{k}(w)}  +   \Vert \theta \Vert^{2/s}_{L^{2}(w)}    \Vert \theta \Vert^{1-\gamma}_{ H^{k}(w)}\Vert \theta \Vert^{1+\gamma}_{ H^{k+\frac{\alpha}{2}}(w)}
\end{eqnarray*}
Since  we have $1+\gamma=2-\ep$, the previous inequality becomes
\begin{eqnarray*}
(II)\lesssim \frac{(2s-1)}{2s} {\mu_3^{\frac{2s}{2s-1}}}\Vert \theta \Vert^{2}_{ H^{k+\frac{\alpha}{2}}(w)} + \frac{1}{2s\mu^{2s}_3 }  \Vert \theta \Vert^{2}_{ H^{k}(w)}  +  \Vert \theta \Vert^{2/s}_{L^{2}(w)}    \Vert \theta \Vert^{\ep}_{ H^{k}(w)}\Vert \theta \Vert^{2-\ep}_{ H^{k+\frac{\alpha}{2}}(w)}.
\end{eqnarray*}
The last term of the previous inequality can be estimated as we did for in inequality \ref{eq: Inegal}. Indeed, it suffices to switch the norms from $\dot H^{k}$  to $\dot H^{k + \frac{\alpha}{2}}$ and to replace $\gamma$ by $\ep$ (this is allowed since the only condition we used throughout these steps was $\gamma>0$) we analogously infer that, for all $\mu_4>0$
\begin{eqnarray*}
&(II)& \lesssim  \left(\frac{(2s-1)}{2s} {\mu_3^{\frac{2s}{2s-1}}}  + \frac{\mu^{1+\ep}_4}{1+\ep} \right) \Vert  \theta \Vert^{2}_{ H^{k+\frac{\alpha}{2}}(w)}   +  \left( \frac{1}{2s\mu^{2s}_3 }  + \frac{1-\ep}{2}  +\frac{\ep}{2\mu_{4}^{1+1/\ep}(1+\ep)^2} \right) \Vert \theta \Vert^{2}_{\dot H^{k}(w)}  \\
 && \ + \  \left( \frac{1-\ep}{2}  +   \frac{\ep(1+2\ep)}{2\mu_{4}^{1+1/\ep}(1+\ep)^2}  \right) \Vert \theta \Vert^{2}_{L^{2}(w)}.
\end{eqnarray*}

\noindent For $(III)$, we have that, for all $\mu_5>0$
$$
(III) \lesssim  \Vert  \theta \Vert_{\dot H^{k+\frac{\alpha}{2}}(w)} \Vert \theta \Vert_{\dot H^{k}(w)} \lesssim  \frac{\mu_5}{2} \Vert  \theta \Vert^{2}_{\dot H^{k+\frac{\alpha}{2}}(w)} +  \frac{1}{2\mu_5}\Vert \theta \Vert^{2}_{\dot H^{k}(w)}.
$$

\noindent Then, using the bound  \ref{L2b},
and,  for all $i\in [0,5]$,  we choose  $\mu_i$   and $\eta$ small enough (recall that $0<1-\gamma =\ep<1$) so that
$$
\frac{1}{2} \partial_{t} \Vert \theta \Vert^{2}_{H^{k}(w)} \leq C \Vert \theta \Vert^{2}_{H^{k}(w)}.
$$
Integrating in time $s \in [0,T]$ and using Gr\"onwall's inequality, we conclude that for all $T<\infty$
$$
\Vert \theta (T) \Vert^{2}_{H^{k}(w)} \leq \Vert \theta_{0} \Vert^{2}_{H^{k}(w)} e^{CT},
$$
where $C$ depends only on $\beta$ and $\Vert \theta_0 \Vert_{L^{\infty}}$.  Moreover, we also obtain
$$
\frac{1}{2} \partial_{t} \Vert \theta \Vert^{2}_{H^{k}(w)} + \int \vert \Lambda^{\frac{\alpha}{2}+k} \theta \vert^{2} \ w dx \leq C(\beta,\Vert \theta_{0} \Vert_{L^{\infty}}) \Vert \theta \Vert^{2}_{ H^{k}(w)}.
$$
Integrating in time $s\in[0,T]$, we obtain
$$
 \Vert \theta(T)\Vert^{2}_{H^{k}(w)} + C \int_{0}^{T} \int \vert \Lambda^{\frac{\alpha}{2}+k} \theta \vert^2 w \ dx \ ds \lesssim \Vert \theta_0\Vert^{2}_{H^{k}(w)} + \int_{0}^{T} \Vert \theta(x,s) \Vert^{2}_{H^{k}(w)} \ ds.
$$
Hence, for all $T<\infty$
$$
 \int_{0}^{T} \int \vert \Lambda^{\frac{\alpha}{2}+k} \theta \vert^2 w \ dx \ ds < \infty.
 $$
 We conclude that for all $T<\infty$, we have $\theta \in L^{\infty}([0,T], H^{k}(w(x) dx)) \cap L^{2}([0,T] , \dot H^{k+\alpha/2}(w(x) dx))$. \\
 
 To conclude the proof of Theorem \ref{sub}, we pass to the weak limit as $R\rightarrow \infty$. To do so, we consider a sequence of solutions $( \theta_{0,m})_{m\in\mathbb{N}^{*}}$.  The strong convergences of the truncated initial data in $H^{2}_{w}$ and in $L^{\infty}$ are straighforward. The {\it{a  priori}} estimates of the previous section allows us to get that the sequence $\theta_m$ is bounded in the space $L^{\infty}([0,T], H^{\alpha/2}(\dw)) \cap L^{2}([0,T], H^{k+\frac{\alpha}{2}}(\dw))$  for all $0<T<\infty$. Then, if $\phi(x,t) \in \mathcal{D}((0,\infty] \times \mathbb R)$, we get that $\phi\theta_m$ is bounded in $L^{2}([0,T], H^{k+\frac{\alpha}{2}}(\dw))$. In order to apply Rellich's theorem we need a bound on $\partial_{t}({\phi}\theta_m)$ or equivatlently on the quantity $\theta_m\partial_{t}\phi  + \phi \partial_{t}\theta_m$. It suffices to focus on $$ \phi \partial_{t}\theta_m=-\phi \partial_{x}(\theta_{m} \mathcal{H}\theta_{m})+\theta_{m} \Lambda \theta_{m}-\Lambda^{\alpha} \theta_{m}.$$
Using the {\it{a  priori}}  $L^{2}([0,T], H^{k+\frac{\alpha}{2}}(\dw))$ bound and the  $L^{\infty}([0,T], L^{\infty}(\dw))$ bound   on $\theta_{m}$ along with the continuity of the Hilbert transform on Sobolev spaces one  easily infer that  $\phi \partial_{t}\theta_m$ is bounded in $L^{2}([0,T], H^{k-\frac{\alpha}{2}}(\dw))$. Therefore, the Rellich compactness theorem (see \cite{PGLR}) allows us to get the existence of a subsquence $\theta_{m_{n}}$ that converges strongly in $L^{2}_{loc} ((0,\infty) \times \mathbb R)$ toward a function $\theta$.  Using that $\theta_{{m_{n}}}$ is a bounded sequence in  $L^{\infty}([0,T], H^{k}(w(x) dx))$ and  on $L^{2}([0,T] , \dot H^{k+\alpha/2}(w(x) dx))$ whose dual spaces are separable Banach spaces, we obtain the *-weak convergence  when $m_n\to+\infty$  of the subsequence $\theta_{m_{n}}$ toward $\theta$  in the spaces  $L^{\infty}([0,T], H^{k}(\dw))$ and  $L^{2}([0,T], H^{k+\frac{\alpha}{2}}(\dw))$. Then, using the strong  $L^{2}_{loc} ((0,\infty) \times \mathbb R)$ one can pass to the weak limit in the nonlinear term, and by standard procedure we also have the convergence of the linear terms in $\mathcal{D}'((0,\infty] \times \mathbb R)$ and we get that the limit is a solution to $\mathcal{T}_{\alpha}$ in the sense of distribution.  \qed

\section{Local existence in $H^{2}(w_{\lambda,\kappa})$ in the case $0<\alpha<1$}

In this section, we prove  Theorem \ref{local}. For $0<\alpha<1$, we consider the following Cauchy problem,
\begin{equation} 
\ (\mathcal{T}_{\alpha}) \ : \\\left\{
\aligned
&\partial_{t}\theta_R+\partial_{x}\theta_R \mathcal{H}\theta_R+ \nu \Lambda^{\alpha}\theta_R = 0 \hspace{2cm} 
\\ \nonumber
& 0<\theta_{0,R}(x)=\theta_{0}(x) \in H^{2}(w_{\lambda,\kappa})  \cap L^{\infty}
\endaligned
\right.
\end{equation}
Since the usual Sobolev embedding $H^{s}({w})(\mathbb R^{n}) \hookrightarrow L^{\infty}(\mathbb R^{n})$ with $s> n/2$ is not anymore true for general $\mathcal{A}_{\infty}$ weights,  one cannot argue as  in \cite{BG} for instance. However, thanks to the pointwise inequality
\begin{equation} \label{eq}
\left(\int \frac{\vert \theta(x)-\theta(y) \vert^{p}}{\vert x-y\vert^{n+s p}} \ dy \right)^{1/p} \leq C \left[ \mathcal{M}\left(\vert \theta -\theta(x) \vert^{p}\right)(x) \right]^{\frac{1-s}{p}}  \times  \left[ \mathcal{M}(\vert \nabla \theta \vert^{q})(x) \right]^{\frac{s}{q}} 
\end{equation}
valid for all $s \in (0,1)$,  $p\in [1,\infty)$ and  $q \geq \frac{pn}{n+p}$, one recover the weighted Gagliardo-Niremberg and the weighted Sobolev inequalities (with $p\ne \infty$). In particular, we can recover  the weighted Sobolev embedding type $\dot H^{2+s/2}(w) \hookrightarrow \dot W^{2,\frac{2}{1-s}}(w)$ for all $s \in (0,1)$. For a proof of \ref{eq}, we refer to the book of Maz'ya \cite{Maz}, p 641.  \\

\noindent We shall study the following truncated equation and get  {\it{a priori}} estimates uniformly in $R>0$
\begin{equation} 
\ (\mathcal{T}_{\alpha,R}) \ : \\\left\{
\aligned
&\partial_{t}\theta_R+\partial_{x}\theta_R \mathcal{H}\theta_R+ \nu \Lambda^{\alpha}\theta_R = 0 \hspace{2cm} 
\\ \nonumber
& 0<\theta_{0,R}(x)=\theta_{0}(x)\phi_{R}(x) \in H^{2}  \cap L^{\infty}
\endaligned
\right.
\end{equation}
Where the function $\phi_{R}(x)$ have been introduced in the end of section \ref{intro}. Then, via (\cite{BG}, Theorem 3.2) we know that there exists a solution $\theta_{R}$ to equation $(\mathcal{T}_{\alpha,R})$ that is sufficiently regular for the product $\theta_{xx} \partial_{t} \theta_{xx}$ to make sense. We shall prove the following energy estimate,  for all $R>0$
\begin{equation} \label{en}
\partial_{t} \Vert \theta_R \Vert_{H^{2}({w})}^{2} + C \Vert \theta_R \Vert^{2}_{\dot H^{2+\frac{\alpha}{2}}}  \lesssim   \Vert \theta_{R} \Vert^{2}_{ H^{2}({w})} +  \Vert \theta_{R} \Vert^{4}_{ H^{2}({w})} + \Vert \theta_R \Vert^{\frac{16}{3}}_{H^{2}(w) }
\end{equation}
We shall omit to write the dependence on $R$ for the sake of readibility. The control of the $L^{2}(w)$ norm is easy to obtain, indeed we have
\begin{eqnarray} \label{ine}
\partial_{t} \int \frac{\theta^{2}}{2}w \ dx + \int \vert \Lambda^{\alpha/2} \theta \vert^{2} w \ dx  &=& - \int \theta^{2} \Lambda \theta w \ dx \ + \int \theta^{2} \mathcal{H}\theta w_x \ dx \\
& -&  \nonumber \int \Lambda^{\alpha/2} \theta [\Lambda^{\alpha/2},w] \theta \ dx 
\end{eqnarray}
 Then, using the first commutator estimate  of lemma \ref{lem31} available for supercritical values of $\alpha$ along with the weighted Sobolev embedding $\dot H^{1/6}({w}) \hookrightarrow L^{3}({w})$ yield the following estimate
\begin{eqnarray*}
\partial_{t} \int \frac{\theta^{2}}{2}w \  dx + \int \vert \Lambda^{\alpha/2} \theta \vert^{2} w \lesssim \Vert \theta \Vert^{2}_{\dot H^{1/6}({w})}\Vert \theta \Vert_{\dot H^{7/6}({w})}+\Vert \theta \Vert^{3}_{\dot H^{1/6}({w})} + \Vert \theta \Vert^{2}_{H^{\alpha/2}(w)} \lesssim \Vert \theta \Vert^{3}_{H^{2}(w)}
\end{eqnarray*}
Note that, in the evolution of the weighted $L^{2}$ norm, we do not need any positivity assumption on the data since the term on the right hand side of \ref{ine} can be  estimated directly via a weighted Sobolev embedding. Indeed, we are dealing with $H^{2}_{w}$ estimates (actually an $H^{1}_{w}$ estimate allows one  to get rid of the positivity assumption). \\

\noindent  The evolution of the homogeneous part is dealt as follows, we have
\begin{eqnarray*}
\frac{1}{2}\partial_{t} \Vert \theta \Vert^{2}_{H^{2}_{\lambda}} + \int  \vert \Lambda^{2+\frac{\alpha}{2}} \theta \vert^{2} \ w \ dx &=&  -\frac{1}{2}\int w_x (\theta_{xx})^{2}  \mathcal{H} \theta  -\frac{3}{2}\int w(\theta_{xx})^{2} \Lambda \theta  + \int w \theta_{xx}  \theta_{x} \Lambda \theta_{x} \\
&&- \int \Lambda^{2+\frac{\alpha}{2}} \theta [\Lambda^{\alpha/2},w] \theta_{xx} \\
&=&I_1+I_2+I_3+I_4
\end{eqnarray*}
In order to estimate $I_1$, we use H\"older's inequality with $\frac{1}{2}=\frac{1}{p}+\frac{1}{q}$ with $p=\frac{2}{1-\alpha}$  and $q=\frac{2}{\alpha}$ and then Young's inequality ($p^{-1}_{1}+p^{-1}_{2}=1$, $p^{-1}_{3}+p^{-1}_{4}=1$). We shall also use the weighted Sobolev embeddings $ H^{\alpha/2}(w) \hookrightarrow L^{\frac{2}{1-\alpha}}(w)$ and $ H^{\frac{1}{2}-\frac{\alpha}{2}}(w) \hookrightarrow L^{\frac{2}{\alpha}}(w)$ that is valid for $\alpha \in (0,1)$ and the fact that $w_{\lambda} \in \mathcal{A}_{\frac{2}{\alpha}}$.  Therefore, we obtain
\begin{eqnarray*}
I_1 &\lesssim& \frac{1}{p_1\ep_1^{p_1} }\Vert \theta \Vert^{p_1}_{\dot H^{2}({w})}  + \frac{\ep_1^{p_2}}{p_2}\Vert \theta_{xx} \Vert^{p_2}_{L^{\frac{2}{1-\alpha}}({w})}\Vert \mathcal{H}\theta \Vert^{p_2}_{L^{ \frac{2}{\alpha}}({w})} \\
&\lesssim&  \frac{1}{p_1\ep_1^{p_1} }\Vert \theta \Vert^{p_1}_{\dot H^{2}({w})}  + \frac{\ep_1^{p_2}}{p_2}\Vert \theta_{xx} \Vert^{p_2}_{H^{\alpha/2}(w)} \Vert \theta \Vert^{p_2}_{H^{\frac{1}{2}-\frac{\alpha}{2}}(w) }  \\
&\lesssim&  \frac{1}{p_1 \ep_1^{p_1} }\Vert \theta \Vert^{p_1}_{\dot H^{2}({w})}  + \frac{\ep_1^{p_2}}{p_3p_2} \Vert \theta \Vert^{p_3p_2}_{\dot H^{2+\frac{\alpha}{2}}({w})} + \frac{\ep_1^{p_2}}{p_4p_2}\Vert \theta \Vert^{p_2p_4}_{H^{2}(w) }
\end{eqnarray*}

 Then, we choose $p_2=1+\eta_2$ and $p_3=1+\eta_3$ where $\eta_2>0$, $\eta_3>0$  are such that  $\eta_2+\eta_3+ \eta_2\eta_3\leq1$ this gives $p_2p_3\leq2$. For instance, one may take $\eta_2=\eta_3=1/3$ hence $p_2=p_3=4/3$ and then $p_1=p_4=4$. Therefore, for all $\ep_1$ we obtain

$$
I_1 \lesssim \frac{1}{4 \ep_1^{4} } \Vert \theta \Vert^{4}_{ H^{2}({w})} + \frac{3\ep_1^{4/3}}{16}\Vert \theta \Vert^{\frac{16}{3}}_{H^{2}(w) }+ \frac{9\ep_1^{4/3}}{16} \Vert \theta \Vert^{2}_{\dot H^{2+\frac{\alpha}{2}}({w})}. \\
$$ 
\noindent To estimate $I_2$ and $I_3$, we follow the same steps as before. We shall keep the same constant to estimate these two terms. \\

\noindent We first observe that $I_2$ is slightly more singular than $I_1$,  however, due to $\Lambda=-\partial_{x}\mathcal{H}$ and the continuity of  $ \mathcal{H}$ on $L^{2/\alpha}(w)$ because  $w \in \mathcal{A}_{\frac{2}{\alpha}}$, following what we did for $I_1$, one gets   

\begin{eqnarray*}
I_2 &\lesssim& \frac{1}{p_1\ep_1^{p_1} }\Vert \theta \Vert^{p_1}_{ H^{2}({w})}  + \frac{\ep_1^{p_2}}{p_2}\Vert \theta_{xx} \Vert^{p_2}_{L^{\frac{2}{1-\alpha}}({w})}\Vert \mathcal{H}\theta_{x} \Vert^{p_2}_{L^{ \frac{2}{\alpha}}({w})} \\
&\lesssim&  \frac{1}{p_1\ep_1^{p_1} }\Vert \theta \Vert^{p_1}_{ H^{2}({w})}  + \frac{\ep_1^{p_2}}{p_2}\Vert \theta_{xx} \Vert^{p_2}_{H^{\alpha/2}(w)} \Vert \theta_{x} \Vert^{p_2}_{H^{\frac{1}{2}-\frac{\alpha}{2}}(w) }  \\
&\lesssim&  \frac{1}{p_1 \ep_1^{p_1} }\Vert \theta \Vert^{p_1}_{ H^{2}({w})}  + \frac{\ep_1^{p_2}}{p_3p_2} \Vert \theta \Vert^{p_3p_2}_{\dot H^{2+\frac{\alpha}{2}}({w})} + \frac{\ep_1^{p_2}}{p_4p_2}\Vert \theta \Vert^{p_2p_4}_{H^{\frac{3}{2}-\frac{\alpha}{2}}(w) }
\end{eqnarray*}
hence, since $H^{2}(w) \hookrightarrow H^{\frac{3}{2}-\frac{\alpha}{2}}(w)$ for $\alpha \in (0,1)$,  we finally get that, for all $\ep_1>0$ 
$$
I_2 \lesssim \frac{1}{4 \ep_1^{4} } \Vert \theta \Vert^{4}_{ H^{2}({w})} + \frac{3\ep_1^{4/3}}{16}\Vert \theta \Vert^{\frac{16}{3}}_{H^{2}(w) }+ \frac{9\ep_1^{4/3}}{16} \Vert \theta \Vert^{2}_{\dot H^{2+\frac{\alpha}{2}}({w})}. \\
$$
To estimate $I_3$, it suffices to observe that this term is as regular as $I_2$. Roughly speaking, these two terms just differ  from an Hilbert transform since we can write
$$
I_3=- \int w \theta_{xx}  \mathcal{H} \theta_{xx}\theta_{x}  \ dx,
$$
therefore,  using the continuity of $\mathcal H$ on $L^{\frac{2}{1-\alpha}}$ because $w \in \mathcal{A}_{\frac{2}{1-\alpha}}$,  we obtain the same estimate as $I_2$, that is, for all $\ep_1>0$ 

$$
I_3 \lesssim \frac{1}{4 \ep_1^{4} } \Vert \theta \Vert^{4}_{ H^{2}({w})} + \frac{3\ep_1^{4/3}}{16}\Vert \theta \Vert^{\frac{16}{3}}_{H^{2}(w) }+ \frac{9\ep_1^{4/3}}{16} \Vert \theta \Vert^{2}_{\dot H^{2+\frac{\alpha}{2}}({w})}. \\
$$

The term $I_4$  can be estimated  thanks to  the first commutator estimate of  lemma \ref{lem31}, then Young's inequality yields
$$
 \vert I_{4} \vert \leq \left\vert\int \sqrt{w}\Lambda^{2+\frac{\alpha}{2}} \theta \ [\Lambda^{\alpha/2},w] \frac{\theta_{xx}}{\sqrt{w}} \ dx \right\vert \lesssim \frac{\ep_2}{2} \Vert \theta \Vert^{2}_{\dot H^{2+\frac{\alpha}{2}}_{w}} + \frac{2}{\ep_2} \Vert \theta \Vert^{2}_{\dot H^{2}_{w}}.
$$
Finally, we have obtained
\begin{equation*} 
\partial_{t} \Vert \theta \Vert_{H^{2}({w})}^{2} +  \Vert \theta \Vert^{2}_{\dot H^{2+\frac{\alpha}{2}}}  \lesssim \left( \frac{\ep_2}{2} + \frac{27\ep_1^{4/3}}{16}\right) \Vert \theta \Vert^{2}_{\dot H^{2+\frac{\alpha}{2}}_{w}} + \frac{2}{\ep_2} \Vert \theta \Vert^{2}_{ H^{2}({w})} + \frac{3}{4 \ep_1^{4} } \Vert \theta \Vert^{4}_{ H^{2}({w})} + \frac{9\ep_1^{4/3}}{16}\Vert \theta \Vert^{\frac{16}{3}}_{H^{2}(w) }
\end{equation*}
Then, choosing $\ep_{1}$ and $\ep_{2}$  \ sufficiently small enough gives the desired {\it{a  priori}} estimate \ref{en}. \\

\noindent  The strong convergence of the initial data as $R\rightarrow \infty$ in $H^{2}$ and $L^{\infty}$ is straighforward. We just focus on the passage to the weak limit in the nonlinearity.  For all $\eta(x,t) \in \mathcal{D}((0,\infty] \times \mathbb R)$, we have that $\eta\theta_k$ is bounded in $L^{2}([0,T], \dot H^{2+\frac{\alpha}{2}}(\dw))$ and then there exists a subsequence $\theta_{k_{n}}$ that converges weakly  to some $\theta$. Moreover,  writting that $ \eta \partial_{t}\theta_k=-\eta\partial_{x}(\theta_{k} \mathcal{H}\theta_{k})+\eta\theta_{k} \Lambda \theta_{k}-\eta\Lambda^{\alpha} \theta_{k}$ we get, by the  {\it{a  priori}} bound in the space  $L^{2}([0,T^{*}], \dot H^{2+\frac{\alpha}{2}}(\dw)) \cap L^{\infty}([0,T^{*}], L^{\infty}(\dw))$ that $\eta \partial_{t}\theta_{k}$ is bounded in $L^{2}([0,T^{*}], \dot H^{1+\frac{\alpha}{2}}(\dw))$. Therefore, one obtains the strong convergence of a subsequence  $\theta_{k_{n'}}$  to a function $\theta$ in $(L^2((0,T^{*}] \times \mathbb R))_{loc}$ which allows to pass to the limit in the nonlinearity.  Then, it is classical to prove that the limit is a solution to the equation in the weak sense.  \\
\qed

\noindent {\bf{Acknowledgment}}:  The author thanks P-G. Lemari\'e-Rieusset  for useful discussions. He also thanks Diego C\'ordoba for useful discussions regarding this model.  He was  supported by the ERC grant Stg-203138-CDSIF and the National Grant MTM2014-59488-P from the Spanish government.

\bibliographystyle{amsplain}

\vspace{1cm}

\begin{quote}
\begin{tabular}{ll}
\vspace{0,2cm}
Omar Lazar & \\
{\small Instituto de Ciencias Matem\'aticas (ICMAT)} & {\small }\\
{\small Consejo Superior de Investigaciones Cient\'ificas} & {\small }\\
{\small C/ Nicolas Cabrera 13-15, 28049 Madrid, Spain} & {\small }\\
{\small Email: omar.lazar@icmat.es} & {\small}
\end{tabular}
\end{quote}
\end{document}